\documentclass[11pt,a4paper,reqno]{amsart}
\usepackage[a4paper,margin=28mm]{geometry}
\usepackage{lmodern,amsmath,amssymb,amsthm,mathtools,microtype}

\usepackage{fontspec}
\newfontface\katakana{KGKatakana.otf}[Path=fonts/,Scale=0.85]
\usepackage{tikz-cd}
\usepackage[numbers]{natbib}
\usepackage[hidelinks]{hyperref}

\newcommand{\kk}{k}
\providecommand{\offf}[1]{\{#1\}}
\hypersetup{pdftitle={Higher cluster-tilting modules for blocks of finite groups},pdfauthor={Tor Kringeland}}

\theoremstyle{definition}
\newtheorem{maintheorem}{Theorem}

\newtheorem{theorem}{Theorem}

\newtheorem{lemma}[theorem]{Lemma}

\DeclareMathOperator{\Hom}{Hom}
\DeclareMathOperator{\End}{End}
\DeclareMathOperator{\Ext}{Ext}
\DeclareMathOperator{\add}{add}

\DeclareMathOperator{\rad}{rad}

\newcommand{\OO}{\mathcal O}

\newcommand{\stHom}{\underline{\Hom}}
\newcommand{\stEnd}{\underline{\End}}
\newcommand{\modu}{\operatorname{mod}}
\newcommand{\dual}{\mathbb D}
\newcommand{\se}{\mathop{\text{\normalfont\katakana セ}}\nolimits}
\newcommand{\tight}{\setlength{\itemsep}{0pt}\setlength{\parsep}{0pt}}

\title[Higher cluster-tilting modules for group blocks]{Higher cluster-tilting modules\\ for blocks of finite groups}
\author{Tor Kringeland}
\date{}

\begin{document}
\begin{abstract}
  We classify the $d$-representation-finite blocks of finite group
  algebras over an algebraically closed field of positive
  characteristic.  The main step is to show that if a block admits a
  $d$-cluster-tilting module with $d\geq3$, then it splits over the
  fraction field of the Witt vectors.  This excludes generalized
  quaternion defect groups of order at least $16$.

\end{abstract}
\maketitle

\section{Introduction}
Which blocks of finite group algebras admit $d$-cluster-tilting
modules, and for which $d$?  For $d\geq 2$ a finitely-generated right
module $M$ over a finite-dimensional algebra $A$ is
$d$-cluster-tilting if
\begin{align*}
  \add M & =\offf{X\in\modu A:\Ext_{A}^{i}(M,X)=0\text{ for }1\leq i<d}  \\
         & =\offf{X\in\modu A:\Ext_{A}^{i}(X,M)=0\text{ for }1\leq i<d}.
\end{align*}
We call $A$ $d$-representation-finite if it has such a
module~\cite{DI}.  Moreover, we write $\se(A)$ for the set of
$d\geq 2$ for which this holds.

Throughout this paper $\kk$ is an algebraically closed field of
characteristic $p>0$, $G$ is a finite group and $B=\kk Gb$ is a block
of $\kk G$ with defect group $D$.  We denote by $l(B)$ the number of
simple modules up to isomorphism.

Known results settle the question for all but one family of defect
groups.  By \citet{EH} and \citet{Bessenrodt} a non-semisimple
$d$-representation-finite block has cyclic or generalized quaternion
defect group, with $d=3$ in the latter case.  The cyclic case follows
from the classification of $d$-representation-finite Brauer tree
algebras in~\cite{DK}, and also from \citet{BHH}.  For $D\cong Q_{8}$
with $l(B)=3$ the block is Morita equivalent to one of two algebras.
\citet{BM} found a $3$-cluster-tilting module over one, and we find
one over the other by a short computation.  The main task of this
paper is to rule out $Q_{2^{n}}$ for $n\geq4$.

\begin{maintheorem}\label{thm:main}
  Let $d\geq 2$.  Then $B$ admits a $d$-cluster-tilting module iff one
  of the following holds:
  \begin{enumerate}\tight
  \item $D=1$;
  \item $D=C_{p}$, $p\geq3$, $l(B)=p-1$ and $d=2p-3$;
  \item $D=C_{7}$, $l(B)\in\offf{3,6}$ and $d=2$; or
  \item $D=Q_{8}$, $l(B)=3$ and $d=3$.
  \end{enumerate}
\end{maintheorem}
For example if $D=C_{7}$ and $l(B)=6$ then (2) holds for $d=11$ and
(3) for $d=2$.  In this case $\se(B)=\offf{2,11}$.

Let us now turn to the remaining quaternion groups, which we exclude
by passing to characteristic zero.  Namely, the block idempotent $b$
lifts to the group algebra $W(\kk)G$ over the Witt vectors, and over
the fraction field $K$ of $W(\kk)$ we get the semisimple algebra
$KGb$.  The algebra $KGb$ is said to be split if all its simple
modules have endomorphism ring $K$.  A character is said to be
$p$-rational if its values lie in $\mathbb Q(\xi)$ for some root of
unity $\xi$ of order prime to $p$.

\begin{maintheorem}\label{thm:split}
  Suppose $B$ has a $d$-cluster-tilting module $M$ for some $d\geq3$.
  Then $KGb$ is split, every irreducible ordinary character of $B$ is
  $p$-rational and $M$ has as many non-isomorphic indecomposable
  summands as there are such characters.
\end{maintheorem}
For $D\cong Q_{2^{n}}$ with $n\geq4$ there are $2^{n-2}-2$ irreducible
characters of $B$ which are not $2$-rational, as shown by
\citet{Olsson}.  Thus we conclude that these blocks have no
$3$-cluster-tilting module.

We recall that Brauer factored the Cartan matrix of a block through
its decomposition matrix.  To prove Theorem~\ref{thm:split} we extend
this factorization from projective modules to rigid ones.  We may
assume that $M$ is basic.  Since $\Ext_{B}^{2}(M,M)=0$ we can lift the
summands of $M$ to the block order $W(\kk)Gb$.  The vanishing of
$\Ext_{B}^{1}(M,M)$ ensures that the $\Hom$-spaces between these
summands retain their dimensions over $K$, and there the matrix of
these dimensions is a Gram matrix of the multiplicities of the simple
$KGb$-modules in the lifts, weighted by the dimensions of the
endomorphism rings of these simple modules.

The same matrix is, up to transpose, the Cartan matrix of
$\End_{B}(M)$.  This algebra has finite global dimension by the higher
Auslander correspondence~\cite{IyamaAC}, so the matrix is unimodular
and its inverse is integral.  Taking traces we find that every weight
is one and that there are as many simple $KGb$-modules as summands of
$M$.  Since $KGb$ is split, the ordinary characters of $B$ take their
values in $K$, and an Eisenstein argument shows that such values are
$p$-rational.  In Section~\ref{sec:lifting} we give the splitting
argument for orders over complete discrete valuation rings, where
neither groups nor symmetry are needed.

\section{Unimodular \texorpdfstring{$\Hom$}{Hom}-matrices force splitting}
\label{sec:lifting}
We compare $\Hom$-spaces between modules before and after lifting them
to an order and extending scalars to the fraction field.  Throughout
this section $\OO$ is a complete discrete valuation ring with
uniformizer $\pi$, residue field $\kk$ (here any field) and fraction
field $K$.  Let $\Lambda$ be an $\OO$-algebra which is free of finite
rank over $\OO$ and whose scalar extension
$K\Lambda=K\otimes_{\OO}\Lambda$ is semisimple.  We put
$A=\Lambda/\pi\Lambda$.  A lift of an $A$-module $X$ is a
$\Lambda$-module $\widetilde X$ which is free of finite rank over
$\OO$, together with an isomorphism
$\widetilde X/\pi\widetilde X\cong X$.  We write $K\widetilde X$ for
$K\otimes_{\OO}\widetilde X$.

Let $M=\bigoplus_{i=1}^{h}X_{i}$ be a basic $A$-module, where the
$X_{i}$ are indecomposable and pairwise non-isomorphic.  We define its
$\Hom$-matrix by
\begin{equation*}
  H_{M}=\bigl(\dim_{\kk}\Hom_{A}(X_{i},X_{j})\bigr)_{1\leq i,j\leq h}.
\end{equation*}
Recall that $M$ is a generator if $A\in\add M$ and that $M$ is rigid
if $\Ext_{A}^{1}(M,M)=0$.  Let $V_{1},\ldots,V_{t}$ be the simple
$K\Lambda$-modules and let
$w_{\alpha}=\dim_{K}\End_{K\Lambda}(V_{\alpha})$ be their weights.
The algebra $K\Lambda$ is said to be split if every weight is one.  If
every $X_{i}$ has a lift $\widetilde X_{i}$ then we define $R$ to be
the $h\times t$ matrix whose entry $R_{i\alpha}$ is the multiplicity
of $V_{\alpha}$ in $K\widetilde X_{i}$.

\begin{theorem}\label{thm:order}
  Let $M=\bigoplus_{i=1}^{h}X_{i}$ be a basic rigid generator of $A$
  such that every $X_{i}$ has a lift.  If $H_{M}$ is unimodular then
  $K\Lambda$ is split, $h=t$ and $\det R=\pm1$.
\end{theorem}
In particular the classes of the $K\widetilde X_{i}$ form a basis of
the Grothendieck group of $K\Lambda$.

\begin{figure}[ht]
  \centering
  \begin{tikzpicture}[x=1.6cm,y=1.1cm,
      d/.style={circle,fill,inner sep=1.4pt}]
    \foreach \x/\t in {0/{\OO},1/{\widetilde P_{+}},2/{\widetilde P_{-}}}
    {\node[d] (u\x) at (\x,1) {};
    \node[anchor=base] at (\x,1.25) {$\t$};}
    \foreach \x/\t in {0/{1},1/{\chi},2/{\operatorname{sgn}}}
    {\node[d] (v\x) at (\x,0) {};
    \node[anchor=base] at (\x,-0.5) {$\t$};}
    \draw (u0)--(v0)--(u1)--(v1)--(u2)--(v2);
  \end{tikzpicture}
  \caption{The lifts of the summands of $M=P_{+}\oplus P_{-}\oplus\kk$
    (top) for $p=3$ and $G=S_{3}$, joined to their ordinary
    constituents (bottom).}
  \label{fig:s3}
\end{figure}
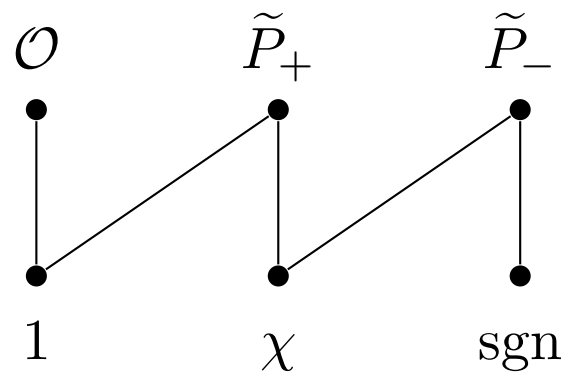

Let us first look at an example with $p=3$, $G=S_{3}$, $\OO=W(\kk)$
and $\pi=p$.  The simple modules of $B=\kk G$ are the trivial module
$\kk$ and the sign module $\kk_{-}$.  Their projective covers $P_{+}$
and $P_{-}$ are uniserial with composition factors $\kk,\kk_{-},\kk$
and $\kk_{-},\kk,\kk_{-}$.  With these uniserial modules one checks
that $M=P_{+}\oplus P_{-}\oplus\kk$ is $3$-cluster-tilting, as in case
(2) of Theorem~\ref{thm:main}.  We lift $P_{+}$ and $P_{-}$ to
projective $\OO G$-modules and $\kk$ to $\OO$ with trivial action.
The ordinary characters $1$, $\operatorname{sgn}$ and $\chi$ (of
degree two) are afforded over $\mathbb Q$, so all weights are one.  On
the $3$-regular elements the character $\chi$ is the sum of the
trivial and sign Brauer characters.  By Brauer reciprocity the three
lifts have characters $1+\chi$, $\operatorname{sgn}+\chi$ and $1$.
Figure~\ref{fig:s3} draws $R$ as a graph.  Each lift is joined to the
constituent below it and, except for $\OO$, to the one to its left, so
$R$ is unitriangular in the order of the figure and the lifts form a
basis of the Grothendieck group, as Theorem~\ref{thm:order} predicts.
Counting common neighbors of the lifts, in the order
$\kk,P_{+},P_{-}$, we get the unimodular matrix
\begin{equation*}
  H_{M}=\begin{pmatrix}1&1&0\\1&2&1\\0&1&2\end{pmatrix}=RR^{\mathsf T}.
\end{equation*}

To see why the unimodularity assumption is needed, take $p=2$,
$G=C_{4}$, $\OO=W(\kk)$ and $\pi=p$.  The rigid generator $M=\kk G$ is
indecomposable, but $H_{M}=(4)$ is not unimodular.  There is no square
root of $-1$ in $K$, since $(1+T)^{2}+1$ is Eisenstein over $\OO$.
Thus $KG\cong K\times K\times K(\sqrt{-1})$, with weights $1$, $1$ and
$2$.  With $R=(1\ 1\ 1)$ the entry of $H_{M}$ counts the constituents
with their weights, $4=1+1+2$.

Now let us prove Theorem~\ref{thm:order}.  The next lemma gives the
two lifting facts needed here and in the proof of
Theorem~\ref{thm:split}.

\begin{lemma}\label{lem:lift}
  Let $X$ and $Y$ be finitely generated $A$-modules.
  \begin{enumerate}\tight
  \item If $\Ext_{A}^{2}(X,X)=0$ then $X$ has a lift.
  \item If $\widetilde X$ and $\widetilde Y$ are lifts of $X$ and $Y$
    and $\Ext_{A}^{1}(X,Y)=0$, then
    \begin{equation*}
      \dim_{\kk}\Hom_{A}(X,Y)=\dim_{K}\Hom_{K\Lambda}(K\widetilde X,K\widetilde Y).
    \end{equation*}
  \end{enumerate}
\end{lemma}
\begin{proof}
  Part (1) holds because the obstructions to lifting $X$ successively
  modulo $\pi^{2},\pi^{3},\ldots$ lie in
  $\Ext_{A}^{2}(X,X)$~\cite[Proposition 1.6]{ADS}.

  For (2) we may reduce a free $\Lambda$-resolution of $\widetilde X$
  modulo $\pi$ to get a free resolution of $X$, since the former
  splits over $\OO$.  Hence
  $\Ext_{\Lambda}^{1}(\widetilde X,Y)\cong\Ext_{A}^{1}(X,Y)=0$.  By
  applying $\Hom_{\Lambda}(\widetilde X,-)$ to the exact sequence
  \begin{equation*}
    0\longrightarrow\widetilde Y\xrightarrow{\ \pi\ }\widetilde Y\longrightarrow Y\longrightarrow0
  \end{equation*}
  and writing down the long exact sequence, we see that multiplication
  by $\pi$ is surjective on
  $\Ext_{\Lambda}^{1}(\widetilde X,\widetilde Y)$.  Since this module
  is finitely generated over $\OO$, it must be zero by Nakayama's
  lemma.  From the same sequence we have $H/\pi H\cong\Hom_{A}(X,Y)$,
  where $H=\Hom_{\Lambda}(\widetilde X,\widetilde Y)$.  Moreover, $H$
  is free over $\OO$ since it is a submodule of
  $\Hom_{\OO}(\widetilde X,\widetilde Y)$, and we have
  $K\otimes_{\OO}H\cong\Hom_{K\Lambda}(K\widetilde X,K\widetilde Y)$,
  so the rank of $H$ gives the required equality.
\end{proof}

\begin{proof}[Proof of Theorem~\ref{thm:order}]
  Write $r_{\alpha}$ for the columns of $R$.  Since $M$ is rigid we
  can compute every entry of $H_{M}$ over $K$ by
  Lemma~\ref{lem:lift}(2).  By Schur's lemma a pair of copies of
  $V_{\alpha}$ contributes $w_{\alpha}$ to a $\Hom$-space, so
  \begin{equation*}
    H_{M}=\sum_{\alpha=1}^{t}w_{\alpha}r_{\alpha}r_{\alpha}^{\mathsf T}.
  \end{equation*}
  Since $A\in\add M$ there is a direct sum $L$ of the
  $\widetilde X_{i}$ with $L/\pi L\cong A$.  Choose $\ell\in L$ that
  maps to $1$ under $L/\pi L\cong A$.  The map $\Lambda\to L$,
  $\lambda\mapsto\ell\lambda$, is an isomorphism modulo $\pi$, hence
  an isomorphism of free $\OO$-modules by Nakayama's lemma.  Thus
  every simple $K\Lambda$-module occurs in some $K\widetilde X_{i}$;
  in other words, every $r_{\alpha}$ is non-zero.

  From the factorization we get that $H_{M}$ is positive semidefinite
  and has rank at most $t$.  As $H_{M}$ is invertible it follows that
  $H_{M}$ is positive definite and $h\leq t$.  Moreover, $H_{M}^{-1}$
  is integral by unimodularity, so it follows that each
  $r_{\alpha}^{\mathsf T}H_{M}^{-1}r_{\alpha}$ is a positive integer.
  Taking traces gives
  \begin{equation*}
    h=\operatorname{tr}(H_{M}^{-1}H_{M})
    =\sum_{\alpha=1}^{t}w_{\alpha}r_{\alpha}^{\mathsf T}H_{M}^{-1}r_{\alpha}
    \geq\sum_{\alpha=1}^{t}w_{\alpha}\geq t\geq h.
  \end{equation*}
  Thus we must have $h=t$ and every $w_{\alpha}=1$, which means that
  $K\Lambda$ is split.  Finally from $\det H_{M}=(\det R)^{2}=\pm1$ we
  get $\det R=\pm1$.
\end{proof}

\begin{proof}[Proof of Theorem~\ref{thm:split}]
  Let $\OO=W(\kk)$ be the ring of Witt vectors over $\kk$, which is a
  complete discrete valuation ring with uniformizer $p$, residue field
  $\kk$ and fraction field $K$ of characteristic zero.  The block
  idempotent $b$ lifts to a central idempotent of $\OO G$, and we also
  denote this idempotent by $b$.  Now the block order $\Lambda=\OO Gb$
  is free of finite rank over $\OO$, and we have
  $\Lambda/p\Lambda\cong B$ and $K\Lambda=KGb$ is semisimple by
  Maschke's theorem.

  Without loss of generality we may assume that $M$ is basic.  We have
  $B\in\add M$ by the defining equalities, since $B$ is projective.
  Since $d\geq3$ we have $\Ext_{B}^{1}(M,M)=\Ext_{B}^{2}(M,M)=0$, and
  every summand of $M$ has a lift by Lemma~\ref{lem:lift}(1).  By the
  higher Auslander correspondence~\cite[Theorem 0.2]{IyamaAC} we know
  that $\End_{B}(M)$ has global dimension at most $d+1$.  As $\kk$ is
  algebraically closed and $M$ is basic, the matrix $H_{M}$ is the
  Cartan matrix of $\End_{B}(M)$ up to transpose.  Since $\End_{B}(M)$
  has finite global dimension, this matrix is unimodular.  From
  Theorem~\ref{thm:order} it follows that $KGb$ is split and $M$ has
  as many indecomposable summands as $KGb$ has simple modules.
  Moreover, the simple $KGb$-modules remain simple over an algebraic
  closure of $K$, since $KGb$ is split.  It follows that they afford
  the irreducible ordinary characters of $B$, all of which take their
  values in $K$.

  It remains to show that a $K$-valued character of $G$ is
  $p$-rational.  To this end let $\xi$ and $\zeta$ be primitive roots
  of unity of orders $|G|_{p'}$ and $|G|_{p}=p^{a}$ in an algebraic
  closure of $K$.  We know that the character values of $G$ lie in
  $\mathbb Q(\xi,\zeta)$.  Since $\kk$ is algebraically closed it
  follows from Hensel's lemma that $\xi\in K$.  In the case $a=0$ we
  are done.  Otherwise we let $\Phi_{p^{a}}$ denote the $p^{a}$th
  cyclotomic polynomial, which has degree $N=(p-1)p^{a-1}$.  Note that
  $\Phi_{p^{a}}(1+T)$ reduces to $T^{N}$ modulo $p$ and has constant
  term $p$, so it is Eisenstein over $\OO$.  It follows that
  $1,\zeta,\ldots,\zeta^{N-1}$ are linearly independent over $K$ and
  that they form a basis of $\mathbb Q(\xi,\zeta)$ over
  $\mathbb Q(\xi)$.  Comparing coefficients in this basis we obtain
  $K\cap\mathbb Q(\xi,\zeta)=\mathbb Q(\xi)$.
\end{proof}

Beware that the hypothesis $d\geq3$ in Theorem~\ref{thm:split} cannot
be weakened to $d\geq2$.  For $p=7$ take the group algebra of the
Frobenius group $C_{7}\rtimes C_{3}$, which is a block with $D=C_{7}$
and $l(B)=3$.  It has a $2$-cluster-tilting module by
Theorem~\ref{thm:main}, but its two characters of degree three take
the values $(-1\pm\sqrt{-7})/2$ at an element of order seven, so they
are not $7$-rational.

\section{Proof of Theorem~\texorpdfstring{\ref{thm:main}}{A}}
\label{sec:classification}
We consider the possible defect groups separately.  If $D=1$ then $B$
is a simple algebra.  In this case every module is projective and
$\add B=\modu B$, whence $B$ is a $d$-cluster-tilting module for all
$d\geq2$.

Now suppose $D\ne1$ and $B$ is $d$-representation-finite.
By~\cite[Theorem 1.1]{EH} every $B$-module has complexity at most one.
On the other hand the maximal complexity of a $B$-module is the
$p$-rank of $D$~\cite[Proposition 2.1]{Bessenrodt}.  Thus $D$ has
$p$-rank one.  The only possibilities are cyclic and generalized
quaternion groups.

Let us first consider $D=C_{p^{a}}$ and put $e=l(B)$.  We have
$e\mid p-1$ by Brauer tree theory, and $B$ is Morita equivalent to a
Brauer tree algebra with $e$ edges and exceptional multiplicity
$(p^{a}-1)/e$~\cite{Dade}.  Such an algebra is
$d$-representation-finite iff either $p^{a}=e+1$, $e\geq2$ and
$d=2e-1$, or $d=2$ and $(p^{a},e)\in\offf{(4,3),(7,3),(7,6)}$,
by~\cite[Corollary 1.3 and Remark 1.4]{DK}.  In the former case we
have $p^{a}=e+1\leq p$ and thus $a=1$, $e=p-1$ and $d=2p-3$, as in
case (2).  The pair $(4,3)$ is impossible because $e\mid p-1$.  This
leaves the other two pairs, which give case (3).

Now suppose that $D$ is generalized quaternion.  Every non-projective
indecomposable $B$-module $X$ satisfies
$\Omega^{4}X\cong X$~\cite[Theorem 5.1]{ES}.  Since $B$ is symmetric,
periodicity and Auslander--Reiten duality give
$\Ext_{B}^{3}(X,X)\cong\stHom_{B}(X,\Omega
  X)\cong\dual\stEnd_{B}(X)\neq0$.  As $B$ is not semisimple, a
$d$-cluster-tilting module has a non-projective summand, so $d\leq3$.
\citet[Section 5.2]{EH} carry out this computation and also exclude
$d=2$.  Hence $d=3$, and by Theorem~\ref{thm:split} every irreducible
ordinary character of $B$ is $2$-rational.  But if $D\cong Q_{2^{n}}$
with $n\geq4$ then some of these characters are not
$2$-rational~\cite[Proposition 4.2]{Olsson}.  This leaves
$D\cong Q_{8}$.

The number of simple modules of a block with defect group $Q_{8}$ is
either one or three.  In the first case $B$ is
nilpotent~\cite[Proposition 3.1]{EEKLS} and hence Morita equivalent to
$\kk Q_{8}$ by \citeauthor{Puig}'s theorem~\cite{Puig}.  Since a
$p$-group algebra has no rigid non-projective module~\cite[Proposition
  2.9]{EKM}, the same is true for $B$.  Thus $B$ cannot have a
$3$-cluster-tilting module.

Finally suppose that $l(B)=3$.  From \citeauthor{Erdmann}'s
classification~\cite[Section IX.5]{Erdmann} it follows that $B$ is
Morita equivalent to one of two algebras, denoted
$Q(3\mathcal A)_{2}^{2}$ and $Q(3\mathcal K)^{2,2,2}$.  Both algebras
are given by quivers with relations over $\mathbb F_{2}$, and we need
to exhibit a $3$-cluster-tilting module for each.  The first algebra
is the basic algebra of the principal block of
$\kk\operatorname{SL}_{2}(5)$, over which \citet[Theorem 1.3]{BM}
found such a module.  The second is the path algebra of the following
quiver
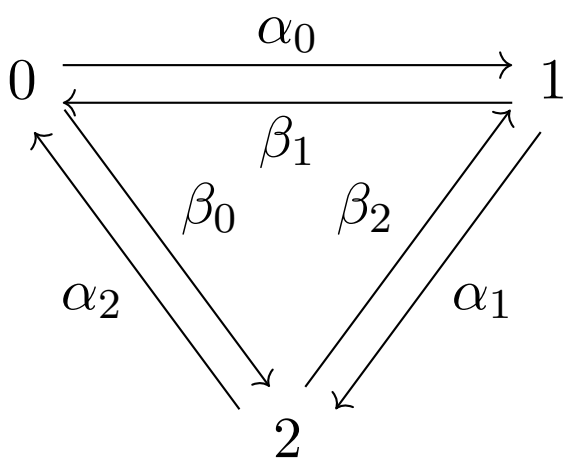
\begin{equation*}
  \begin{tikzcd}[column sep=3.2em,row sep=4.8em,
      every arrow/.append style={shift left=0.75ex},
      every label/.append style={font=\normalsize}]
    0 \arrow[rr,"\alpha_{0}"] \arrow[dr,"\beta_{0}"]
    & & 1 \arrow[ll,"\beta_{1}"] \arrow[dl,"\alpha_{1}"] \\
    & 2 \arrow[ul,"\alpha_{2}"] \arrow[ur,"\beta_{2}"] &
  \end{tikzcd}
\end{equation*}
modulo the relations
\begin{equation*}
  \alpha_{i}\alpha_{i+1}=\beta_{i}\alpha_{i-1}\beta_{i},\qquad
  \beta_{i}\beta_{i-1}=\alpha_{i}\beta_{i+1}\alpha_{i},\qquad
  \alpha_{i}\beta_{i+1}\beta_{i}=0,
\end{equation*}
where the indices are taken modulo $3$ and paths are composed from
left to right.  We denote this algebra by $A$ and the trivial path at
$i$ by $e_{i}$, and take
\begin{equation*}
  M=A\oplus e_{0}A/\alpha_{0}\beta_{1}A\oplus
  e_{0}A/\beta_{0}\beta_{2}A\oplus e_{1}A/\beta_{1}\alpha_{0}A\oplus
  \rad e_{2}A.
\end{equation*}

As $A$ is self-injective and is a summand of $M$, we have that $M$ is
a generator--cogenerator.  By \citeauthor{IyamaAC}'s
criterion~\cite[Theorem 5.1(3)]{IyamaAC} the module $M$ is
$3$-cluster-tilting iff $\Ext_{A}^{1}(M,M)=\Ext_{A}^{2}(M,M)=0$ and
$\End_{A}(M)$ has global dimension at most four.

We let Auslander~\cite{Auslander} do the legwork \textit{in silico}
over $\mathbb F_{2}$.  We find that the two $\Ext$-groups are zero.
We also use minimal right $\add M$-approximations to construct a
minimal projective resolution for each simple $\End_{A}(M)$-module.
Each of these resolutions has length at most four, and four of the
seven have length exactly four.  Moreover we were able to verify both
of these results independently in QPA~\cite{QPA}.  Now since $\Hom$-
and $\Ext$-groups commute with field extensions and global dimension
is preserved under scalar extension from the perfect field
$\mathbb F_{2}$~\cite[Lemma 1.1(2)]{BM}, this computation gives the
same conclusions over $\kk$.  So we conclude that $B$ has a
$3$-cluster-tilting module, as asserted in (4).

\bibliographystyle{plainnat}
{\small\raggedright
  \bibliography{references}}
\end{document}